\documentstyle{article}

\input tcilatex              
\input xypic          

\makeatletter
\def\@mybiblabel#1{#1.}
\renewenvironment{thebibliography}[1]
     {\section*{\refname
        \@mkboth{\MakeUppercase\refname}{\MakeUppercase\refname}}%
      \list{\@mybiblabel{\@arabic\c@enumiv}}%
           {\settowidth\labelwidth{\@biblabel{#1}}%
            \leftmargin\labelwidth
            \advance\leftmargin\labelsep
            \usecounter{enumiv}%
            \let\p@enumiv\@empty
            \renewcommand\theenumiv{\@arabic\c@enumiv}}%
      \sloppy
      \clubpenalty4000
      \@clubpenalty \clubpenalty
      \widowpenalty4000%
      \sfcode`\.\@m}
     {\def\@noitemerr
       {\@latex@warning{Empty `thebibliography' environment}}%
      \endlist}
\makeatother

\def\Hom{{\mathit Hom}}

\def\Hoch{{\mathit Hoch}}

\newcommand{\ot}{\otimes}

\newcommand{\To}{\Rightarrow}
\newcommand{\lo}{\leftarrow}

\newcommand{\ii}{\'{\i}}

\newcommand{\wh}{\widehat}
\newcommand{\wt}{\widetilde}
\newcommand{\dem}{{\bf Proof: }}
\newcommand{\rem}{\noindent{\bf Remark: }}
\newcommand{\rems}{\noindent{\bf Remarks: }}
\newcommand{\ex}{\noindent{\bf Example: }}
\newcommand{\exs}{\noindent{\bf Examples: }}

\newtheorem{teo}{Theorem}[section]
\newtheorem{defi}[teo]{Definition}
\newtheorem{lem}[teo]{Lemma}
\newtheorem{prop}[teo]{Proposition}
\newtheorem{coro}[teo]{Corollary}

  \def\C{{\cal C}} \def\D{{\cal D}}
\def\E{{\cal E}}   \def\H{{\cal H}}
   
\def\M{{\cal M}}   

\font \Bbb=msbm10 scaled \magstep1

 \def\ZZ{\hbox{\Bbb Z}}
\def\NN{\hbox{\Bbb N}} 

\font \bb=msbm10 scaled \magstep0

\def\zz{\hbox{\bb Z}} \def\nn{\hbox{\bb N}}

\font\tenfrak=eufm10 scaled \magstep1
\font\sevenfrak=eufm7 scaled \magstep1
\newfam\frakfam
\textfont\frakfam=\tenfrak
\scriptfont\frakfam\sevenfrak
\def\frak{\fam\frakfam\tenfrak}

\def\Hoch{{\mathit Hoch}}    \def\Hom{{\mathit Hom}}
\def\Chain{{\mathit Chain}}  \def\Vec{{}_k{\mathit Vect}}

\title{On the derived invariance of cohomology theories
for coalgebras}

\author{Marco A. Farinati${}^1$}

\date{}
\begin{document}\maketitle

{\footnotetext[1]
{Dto. de Matem\'atica FCEyN, Universidad de Buenos Aires,
Cdad. Universitaria Pab I. (1428) Buenos Aires - Argentina.\\
e-mail:  mfarinat@dm.uba.ar\\
Research partially supported
CONICET, Fundaci\'on Antorchas and UBACYT TW169}
}

\begin{abstract}
We study the derived invariance of the cohomology
theories  $\Hoch^*$, $H^*$ and $HC^*$ 
associated to coalgebras
over a field. We prove a theorem characterizing
derived equivalences. As particular cases, it
describes the
two following situations: 1) $f:C\to D$ a quasi-isomorphism
of differential graded coalgebras, 2) the existence
of a ``cotilting" bicomodule ${}_CT_D$. In these two cases
we construct a derived-Morita equivalence context, and
consequently we obtain isomorphisms $\Hoch^*(C)\cong\Hoch^*(D)$ and 
$H^*(C)\cong H^*(D)$. Moreover, when we have a coassociative 
map inducing an isomorphism $H^*(C)\cong H^*(D)$ (for example 
when there is a quasi-isomorphism $f:C\to D$), we prove that
$HC^*(C)\cong HC^*(D)$.
\end{abstract}

\tableofcontents

\section*{Introduction}

This work is devoted to the study of cohomology theories
for coalgebras and their relations with the derived 
categories associated to coalgebras.
These cohomology theories are called $\Hoch^*$, $H^*$ 
(defined by Y. Doi in \cite{doi}) and 
$HC^*$ (defined by A. Solotar and myself in \cite{hc}).
They are important categorical invariants,
for instance $H^*$ and $\Hoch$ measure coseparability, and
$H^*$ measures extensions, and hence it is a useful object when studying 
classification problems. Also in some examples they
have a nice geometric interpretation:
in \cite{hc} we prove that the (topological) coalgebra $\D(X)$ 
of distributions over a real compact smooth manifold $X$ has
the $n$-currents (i.e. the continuous dual of the $n$-differential
forms) as the $n$-th $\Hoch$ group; we also show in \cite{hc},
that the cyclic cohomology of $\D(X)$ 
can be computed in terms of the De Rahm cohomology of 
the manifold $X$.

The notion of  coalgebra in a monoidal category $(\C,\ot)$  
coincides with 
the concept of algebra in the category
$(\C^{op},\ot)$, so under this duality, one can translate
``algebraic" statements in order to get ``coalgebraic"
statements and viceversa. But, even if the
notion of (abelian and) monoidal category is
selfdual, in practice, a category $\C$ is usually far from being 
equivalent to $\C^{op}$. If the category
$\C$ is fixed, simple dualization of
``algebraic" proofs needs not lead to ``coalgebraic" proofs,
because one wants theorems to hold in $\C$ and not in $\C^{op}$.
The most familiar example of this situation is $\Vec$, the category
of vector spaces over some fixed ground field $k$, which is
not equivalent to its opposite category (it satisfies 
Grothendieck's axiom AB5 and it is well known that
if an abelian category satisfies AB5 and AB5* then
every object is the zero object).

One of the guiding lines of this work is the
following informal idea: ``the derived category associated
to a coalgebra keeps the homological relevant information, so if
two coalgebras have equivalent derived categories, they should
have isomorphic cohomology". A precise statement of this idea is
Theorem \ref{morita}, together with our two main examples:
quasi-isomorphisms (Proposition \ref{propquis})
and cotilting bicomodules (Theorem \ref{teocotilting}).
The language of differential graded objects appears
naturally in the context of derived categories, so we treat with
differential graded coalgebras, even also when we are looking
for results on usual coalgebras.

We remark that, for a fixed field $k$,
the category of differential $\ZZ$-graded $k$-vector spaces
is isomorphic (as monoidal category) to a category of comodules
over a Hopf algebra (see the first part of section 1), an so
it is not equivalent to its opposite category ( same 
reasons that $\Vec\not\cong \Vec^{op}$).

\medskip

The organization of this work is the following:

In Section 1 we introduce the notion of a differential graded
coalgebra $C$ and the category of differential graded
comodules, denoted by $\Chain(C)$. This category has a natural
notion of homotopy, so we define the homotopy  category $\H(C)$ and
its corresponding localization by quasi-isomorphisms $\D(C)$.
We define the notion of {\em closed} object, which plays the
role of injective resolution. The main result of this
section is Theorem \ref{teo1}, which states that,
when $C$ is positively graded, and we restrict to the category
$\Chain^+(C)$ of left-bounded differential
graded comodules,  then the standard
resolution $C(M)$ is a closed object, quasi-isomorphic to $M$.

Section 2 deals with  derived equivalences. We begin with a general
characterization of equivalences given by Proposition
\ref{propb}
and  we prove Theorem \ref{cokeller}. It characterizes
equivalences induced by derived cotensor products.

If we specialize Theorem \ref{cokeller} to the case of usual coalgebras
and comodules, the hypothesis lead us
to the notion of cotilting bicomodule. In Section
3 we define the notion of cotilting bicomodule
and prove  that they induce equivalences of derived categories
(Theorem \ref{Tbox}), and moreover, that they induce 
derived Morita contexts (Theorem \ref{teocotilting}).

On the other side, concerning non trivial differential
graded structures, we prove in Section 4 that if $f:C\to D$
is a quasi-isomorphism of positively graded differential
coalgebras, then their derived categories are equivalent, and
they fit into a derived Morita context.

Finally Section 5 is devoted to the extension of the definition
of $\Hoch^*$, $H^*$ and $HC^*$ to the differential graded case
extending some properties of the non differential graded case 
(Proposition \ref{sbi}). We prove Theorem \ref{morita}.
This result implies the invariance of $\Hoch^*$ and $H^*$
under quasi-isomorphisms and cotilting equivalences, and
Proposition \ref{propf}, which implies the invariance
of $HC^*$ under quasi-isomorphisms.

We mention that the dualization of the notion of tilting 
modules was introduced in this work with the aim
of understanding all possible (or at least most of the) 
derived equivalences for usual coalgebras. We wanted to answer
the following question: given two $k$-coalgebras $C$ and $D$ such 
that $\D^+(C)\cong \D^+(D)$ as triangulated categories, is it true that
there exists an equivalence (may be another) given by a derived
cotensor product? One may view 
Proposition \ref{propb} and
Theorem \ref{teocotilting} as partial results in this direction.

On the other hand, the example of quasi-isomorphisms of differential
graded coalgebras has a computational motivation: if $C$ is a
subcoalgebra of $C'$, one can try to find a differential graded
coalgebra of the form $C'\ot \Lambda$ where $\Lambda_0=k$ and 
a quasi-isomorphism $C\to C'\ot \Lambda$, in other words,
find a positively graded coalgebra $\Lambda$ and a differential 
$d$ on $C'\ot \Lambda$ such that the cohomology of
$(C'\ot \Lambda,d)$ is zero for positive degrees and the kernel
of $d_0$ is $C$. In this situation we can interpret 
$(C'\ot \Lambda,d)$ as a ``model" for $C$, and Theorem
\ref{morita} and Proposition \ref{propf} tell us that, in order
to compute $H^*$, $\Hoch^*$ or $HC^*$, one
can replace the coalgebra $C$ by its differential graded model.

\medskip

Throughout this paper $k$ will denote a field of
arbitrary characteristic, all unandorned tensor products will be taken 
over $k$.

\bigskip

\noindent{\bf Acknowledgements:} I  wish to thank my Ph.D 
advisor Andrea Solotar for her continuous encouragement
and patience, most of the results of this work are 
part of my Ph.D. thesis \cite{tesis}.
I also wish to thank Bernhard Keller for 
general comments and corrections of preliminary versions
of my thesis.

\section[A triangulated category associated to a D.G.
coalgebra]{A triangulated category associated to a differential
graded coalgebra}

\subsection[D.G. coalgebras and D.G. comodules]
{Differential graded coalgebras and differential 
graded comodules}

The category of differential $\ZZ$-graded $k$-vector spaces
is a monoidal category with $\ot_k$ as tensor product, so the notion
of coalgebra in this category makes sense. By definition, a
{\bf differential graded $k$-coalgebra} is a coalgebra
in the monoidal category of differential graded $k$-vector spaces,
namely, it is the data $(C,\Delta)$ where
$C=(\oplus_{n\in\zz}C_n,d)$ is a differential graded vector
space, $d_n:C_n\to C_{n+1}$, and
$$\Delta:C\to C\ot C$$
is a coassociative map (in the category of differential graded
$k$-vector spaces)
admitting a counit,
i.e.
$\Delta(C_n)\subset (C\ot C)_n=\oplus_{p+q=n}C_p\ot C_q$
and $\Delta d_C =d_{C\ot C}\Delta$ where
$d_{C\ot C}$ is defined on homogeneous elements by:
$$d(c_1\ot c_2):=d(c_1)\ot c_2+(-1)^{|c_1|}c_1\ot d(c_2)$$

It is well-known that, for any group $G$, the category of 
$G$-graded vector spaces is isomorphic to the category of 
$k[G]$-comodules as monoidal categories (with the diagonal
structure in the tensor product), where $k[G]$ denotes the
group-algebra. Also, the category of vector spaces together
with an square zero endomorphism is isomorphic to the category 
of $k[d]/\langle d^2\rangle$-modules. From these two remarks,
the following is a natural way of describing the category
of differential graded $k$-vector spaces
as the category of comodules over a Hopf algebra:

Let $H$ denote the $k$-algebra generated by $X$, $X^{-1}$, $d$
with relations
$$XX^{-1}=X^{-1}X=1$$
$$d^2=0$$
$$Xd+dX=0$$
Then, $H$ admits a Hopf algebra structure determined by
$$\Delta(X)=X\ot X$$
$$\Delta(d)=1\ot d+d\ot X$$
The counit is given by evaluation of $X$ in $1$ and $d$ in $0$,
the antipode by $S(X)=X^{-1}$ and
$S(d)=X^{-1}d=-dX^{-1}$.

Now the definition of differential
graded coalgebra can be given in terms of $H$, because
a coalgebra in ${}^H\M$ is the same
as an $H$-comodule coalgebra. The natural corepresentations
of a differential graded coalgebra $C$ are the so-called
differential graded comodules (i.e. $C$-comodules which are
differential graded vector spaces and the structure map is a map
of differential graded vector spaces). This category can
be identified
with ${}^C\left({}^H\M\right)\cong {}^{C\#H}\M$, where
$C\# H$ is the smash product of $C$ with $H$, and the equivalence
${}^C\left({}^H\M\right)\cong {}^{C\#H}\M$ preserves the underlying
vector space and it is the identity on the arrows.

As a corollary of this simple remark, we have that the category 
of differential graded left comodules over a differential graded 
coalgebra $C$, which will be denoted by $\Chain(C)$, is a 
comodule-category (over a field). Therefore it is a Grothendieck 
category, in particular it has arbitrary sums and products, 
projective limits, filtered inductive limits are exact, and 
each object is the union of its finite dimensional subobjects.

\subsection{The categories $\H(C)$ and $\D(C)$}

The category $\Chain(C)$  
has a natural notion of homotopy.
We will say that two maps $f,g:M\to N$ in $\Hom_{Ch(C)}(M,N)$ are 
{\bf homotopic} and we will write $f\sim g$, in case that there 
exists a graded $C$-colinear map $h:M\to N[1]$ such that 
$f-g=hd_M+d_Nh$. As usual, this is an equivalence
relation in $\Hom_{Ch(C)}(M,N)$, compatible
with addition and composition, so we can define the quotient
category $\H(C)$. By definition, the objects
of $\H(C)$ are the same as the objects of $\Chain(C)$, and the
maps are homotopy clases of maps, i.e. $\Hom_{\H(C)}(M,N)=
\Hom_{Ch(C)}(M,N)/\sim$.

Let's denote by $C^*$ the graded dual of $C$; $C^*$ is 
a differential graded algebra. 
The category $\Chain(C)$ embeds in $\Chain(C^*)$, the category
of differential graded (right) $C^*$-modules.
The notion of a Cone of a morphism can be defined as the
Cone in $\Chain(C^*)$, namely, if $f \in \Hom_{Ch(C)}(M,N)$,
$Co(f)\in\Chain(C)$ is defined  by:
$$Co(f)_n:=M[-1]_n\oplus N_n$$
$$d_{Co(f)}:=d_{M[-1]}+f+d_{N}$$
The cone of 
the identity
map is not the zero object, but homotopically equivalent to it, so 
$\Chain(C)$ is never a triangulated category (with cones as
distiguished triangles). The additive category $\H(C)$ is triangulated, 
taking as distinguished triangles the uples isomorphic to 
$$(M,N,Co(f),f:M\to N,\pi:Co(f)\to M[-1],
f[-1]:M[-1]\to N[-1])$$
The embedding $\H(C)\to\H(C^*)$ is an embedding of
triangulated categories.

The `derived' category associated to $C$, denoted by $\D(C)$, is
the localization of $\H(C)$ by the class of quasi-isomorphisms
(i.e. morphisms in $\H(C)$ inducing an isomorphism in
homology). The category $\D(C)$ is also a triangulated category
with the triangulated structure induced by the localization
functor and the structure of $\H(C)$.

As usual, we will denote $\Chain^{+,-,b}(C)$, $\H^{+,-,b}(C)$ and 
$\D^{+,-,b}(C)$ the full subcategories of $\Chain(C)$, $\H(C)$ and
$\D(C)$ consisting on objects which are isomorphic to left (resp. 
right, both sided) bounded complexes.

We will use the following  Lemma of triangulated categories whose 
proof is obtained adapting that of \cite{rusos} (Lemma 13, chap. 
IV, \S 1) to the differential graded case.
The only thing to check is that every morphism defined in the proof
of their Lemma belongs to the category $\Chain(C)$.
\begin{lem}\label{lemaruso}
Let $C$ be a differential graded coalgebra and
$$0\to X\to Y\to Z\to 0$$
a short exact sequence in $\Chain(C)$ which splits as a
sequence of graded $C$-comodules. Then the sequence
fits into a triangle $X\to Y\to Z\to X[-1]\to\dots$
\end{lem}

\subsection{Characterization of $\D^+(C)$}

From now on, $C$ will denote a differential graded $k$-coalgebra.
\begin{defi}
We will say that an object $X\in \Chain(C)$ is {\bf closed} if,
for all $M\in\Chain(C)$, the localization
funtor gives an isomorphism
$$\Hom_{\H(C)}(M,X)\cong \Hom_{D(C)}(M,X)$$
\end{defi}
As a consequence of the definition, being closed is stable
under products, traslation, and (by a five Lemma argument) triangles, 
namely if two of the three objects in a triangle are closed objects, 
then so is the third one.

We will see later that there exist enough closed objects,
at least when $C$ is positively graded and we restrict ourselves
to the subcategory $\Chain^+(C)$.
We begin with an adjunction:
\begin{prop}\label{adjuncion}
Let $V \in \Chain(k)$, $M \in\Chain(C)$, and consider
$C\ot V \in\Chain(C)$ with the structure of $C$ comodule given by
$C$ and the standard graded differential structure of the tensor 
product, then
$$\Hom_{\H(C)}(M,C\ot V)\cong \Hom_{\H(k)}(M,V)=$$
$$=\Hom_{\D(k)}(M,V)\cong \Hom_{\D(C)}(M,C\ot V)$$
\end{prop}
\dem The well-known adjunction in the standard comodule context can 
be carried out in the differential graded case, defining maps
$$\Hom_{Ch(C)}(M,C\ot V)\to \Hom_{Ch(k)}(M,V)$$
$$f\mapsto (\epsilon\ot id_V)\circ f$$
and
$$\Hom_{Ch(k)}(M, V)\to \Hom_{Ch(C)}(M,C\ot V)$$
$$g\mapsto (id\ot g)\circ \rho_M$$
where $\rho_M:M\to C\ot M$ denotes the structure map of $M$.

If $f,f':M\to C\ot V$ are homotopic maps
by means of an homotopy $h$, then it is easily checked
that $(\epsilon\ot id)h$ is an homotopy between
$(\epsilon\ot id)f$ and $(\epsilon\ot id)f'$ (one uses that $\epsilon$
is necessarily a graded morphism and that $\epsilon \circ d=0$).
Conversely, if $g,g':M\to V$ are two homotopic maps, let's call
again $h$ the homotopy between $g$ and $g'$, then
$(id_C\ot h)\circ \rho_M$ gives the desired homotopy between
$(id_C\ot h)\circ g$ and $(id_C\ot h)\circ g'$, and
the same formula as above gives a well-defined isomorphism
$$\Hom_{\H(C)}(M,C\ot V)\cong \Hom_{\H(k)}(M,V)$$
Now if $f:M\to M'$ is a quasi-isomorphism of differential
graded $C$-comodules, then obviously the forgetful functor gives
a quasi-isomorphism, and considering the other variable of
the $\Hom$, if $\phi:V\to V'$ is a quasi-isomorphism of differential
graded vector spaces, then by the K\"unneth formula
$id_C\ot \phi:C\ot V\to C\ot V'$ is also a quasi-isomorphism.
As a consequence, the adjointness remains valid after localization
giving 
$$\Hom_{\D(C)}(M,C\ot V)\cong \Hom_{\D(k)}(M,V)$$
Finally, we remark that every quasi-isomorphism in $\Chain(k)$ 
is an homotopy equivalence, then $\D(k)=\H(k)$.

\bigskip

As a consequence of the above proposition, for every differential
graded vector space $V$, the object $C\ot V\in\Chain(C)$ is closed.

We will state a theorem which characterizes the category $\D^+(C)$
for positively graded differential coalgebras:

\begin{teo}\label{teo1}
Let $C$ be a positively graded differential coalgebra and
let $M$ be an object in $\Chain^+(C)$. Then, there is a 
functorial way of assigning to $M$ a closed object $C(M)$ and
a quasi-isomorphism $M\to C(M)$; the object $C(M)$ also belongs to
$\Chain^+(C)$.
\end{teo}
As a consequence of this theorem, the category $\D^+(C)$ can be
described as the full subcategory of $\H^+(C)$ consisting
of closed objects. The proof of the theorem is achieved by
exposing an explicit standard resolution, and proving that
this resolution is a closed object. In order to do that, we need 
two Lemmas. We begin by defining the standard resolution:

\medskip

Let $M\in \Chain(C)$ where $C$ is any differential graded coalgebra,
we define the object $C(M)\in\Chain(C)$ by
$C(M):=\oplus_{n>0}C^{\ot n}\ot M$ with the following graduation
and differential:
$$C(M)_p=\bigoplus_{i_1+\dots+i_r+j+r-1=p}C_{i_1}\ot\dots C_{i_r}\ot M_j$$
and, for $(c_{i_1},\dots,c_{i_r},m)\in
C_{i_1}\ot\dots C_{i_r}\ot M$,
$$\partial(c_{i_1},\dots,c_{i_r},m)=
(-1)^r\sum_{k=1}^{r+1}(-1)^{|c_{i_1}|+\dots+|c_{i_{k-1}}|}(c_{i_1},\dots,
d(c_{i_k}),\dots,c_{i_r},m)+$$
$$+\sum_{k=1}^{r+1}(-1)^{k+1}(
c_{i_1},\dots, \Delta(c_{i_k}),\dots,c_{i_r},m)$$
where we use the convention $c_{i_{r+1}}=m$, and
$\Delta(m)=\rho_M(m)$. In an abridged way, we write
$$\partial(c_{i_1},\dots,c_{i_r},m)=
d(c_{i_1},\dots,c_{i_r},m)
+b'(c_{i_1},\dots,c_{i_r},m)$$
\begin{lem}
With the above notations, $\partial^2=0$, and the structure morphism
$\rho_M:M\to C\ot M\subseteq C(M)$ is a quasi-isomorphism in
$\Chain(C)$.
\end{lem}

\dem Let us define the extended complex  
$\wh{C}(M)=(\oplus_{n\geq 0}C^{\ot n}\ot M,b',d)$ with differential 
given by the same formulae as in $C(M)$ and the $k$-linear map
$h:C(M)\to C(M)$ by
$$h(c_1,\dots ,c_r,m)=\epsilon(c_1)(c_2,\dots,c_r,m)
\ \ \hbox{ for $n\geq 1$}$$
and $h(m)=0$ for $m\in M$. 

It is well known that $b'^2=0$, $d^2=0$
and that $b'h+hb'=id$, then it is sufficient to see that
$b'd+db'= 0 = dh+hd$. Since $\wh{C}(M)$ is the cone of 
the map $\rho_M:M\to C(M)$, the equality 
$b'd+db'= 0$ says that
$\partial$ is in fact a square zero operator,
and the equality $0 = dh+hd$ 
means that $\partial h+h\partial = id$,
which implies that $\wh{C}(M)$ is acyclic, or equivalently that
$\rho_M:M\to C(M)$ is a quasi-isomorphism.

We omit the computation of  $b'd+db'$ and $dh+hd$, to see that they 
equal zero is tedious but straightforward.

\bigskip

It is clear that the assignment $M\mapsto C(M)$ is functorial, and 
by the previous Lemma, $M$ is quasi-isomorphic to $(C(M), \partial)$. 
We do not 
know if $C(M)$ is a closed object in general, but we succeeded to 
prove it (as it was enounced in Theorem \ref{teo1})
for $C$ positively graded and $M$ bounded below, i.e.
for $M$ and object in $\Chain^+(C)$. The proof of this result
follows from the following Lemma:

\begin{lem}\label{tri}
Let $\{M_{n}\}_{n\in\nn}$ be an inverse system in $\Chain(C)$ with 
morphisms $p_{n+1}:M_{n+1}\to M_n$ which are, as morphisms of graded 
$C$-comodules, split epimorphisms. Then, the short exact sequence
$$0\to \stackunder{\lo n}{\lim}M_n \to
\prod_{n\in \nn}M_n\to \prod_{n\in \nn}M_n\to 0$$
fits into a triangle in the category $\H(C)$
where the products (and the inverse limit) are taken in the category
$\Chain(C)$, the map between the products is
$$\{m_n\}_{n\in\nn}\mapsto\{m_n-p_{n+1}(m_{n+1})\}_{n\in\nn}$$
\end{lem}
\dem let us call $s_n:M_n\to M_{n+1}$ the splittings of  the maps $p_n$,
and define $s:\prod_n M_n\to \prod_n M_n$ by
$$s(m_0,m_1,\dots,m_n,\dots)=\ \ $$
$$=
(m_0,0,-s_1(m_1),-s_2(s_1(m_1))-s_2(m_2),
-s_3(s_2(s_1(m_1)))-s_3(s_2(m_2))-s_3(m_3),\dots)$$

This proves that the sequence
$$0\to \stackunder{\lo n}{\lim}M_n \to
\prod_{n\in \nn}M_n\to \prod_{n\in \nn}M_n\to 0$$
splits as a sequence of graded $C$-comodules. Now use Lemma
\ref{lemaruso} and the proof is completed.

\bigskip

Now consider, given $M\in\Chain^+(C)$, the system
$C_n(M):=(\oplus_{i=1}^nC^{\ot i}\ot M,b',d)$ with differential
induced by the projection $C(M)\to C_n(M)$. For every $n\in\NN$
there is a graded $C$-split short exact sequence
$$0\to C^{\ot n+1}\ot M \to C_{n+1}(M)\to C_n(M)\to 0$$
which implies two facts:
\begin{itemize}
\item When $n=1$, $C_1(M)=C\ot M$ is a closed object; for $n>1$,
$C^{\ot n}\ot M=C\ot (C^{\ot n-1}\ot M)$ is also a closed object,
and hence, the above exact sequence implies inductively that $C_n(M)$
is closed for all $n\in\NN$.
\item The maps of the system $C_{n+1}(M)\to C_{n}(M)$,
viewed as morphisms of graded $C$-comodules, are
split surjections, and the hypothesis of Lemma \ref{tri}
are fulfilled.
\end{itemize}
As a consequence, we can conclude that 
$C(M)=\stackunder{\lo n}{\lim}C_n(M)$ is a closed object, and the
proof of  Theorem \ref{teo1} is complete.

\bigskip

We remark that if $C$ has nonzero components in negative degrees, or 
if $M$ is a left  unbounded complex (i.e. has infinitely many nonzero
components of negative degree), then
$C(M)\neq \stackunder{\lo n}{\lim}C_n(M)$. In general
$\stackunder{\lo n}{\lim}C_n(M)
= Tot^{\prod}(C^{\ot *}\ot M,b',d)$,
which is closed (by the same arguments) but 
unfortunately is not necessarily
quasi-isomorphic to $M$.

\bigskip

Theorem \ref{teo1} gives also a way of defining derived functors
for each functor between the homotopy categories. We define in 
this way the derived cotensor product, 
noticing that, given ${}_DX_C\in \Chain(D\ot C^{op})$
(where $C$ and $D$ are differential graded coalgebras),
$X\Box_C-:\Chain(C)\to \Chain(D)$ preserves homotopies,
then it defines a functor $X\Box_C-:\H(C)\to\H(D)$. In order
to have a derived cotensor product between $\D^+(C)$ and
$\D^+(D)$ we need $C$ and $D$ to be positively graded, and
$X\in \Chain^+(D\ot C^{op})$. 

Also, following the proof of Theorem \ref{teo1}, we have 
the following corollary:

\begin{coro}
Let $C$ be a positively graded differential coalgebra, $X\in\Chain^+(C)$
a closed object and $V\in\Chain^+(k)$, then $X\ot V$ is also a closed
object. 
\end{coro}
\dem we know after Theorem \ref{teo1} that, for any
object $M\in\Chain^+(C)$, $C(M)$ is closed. Since $\rho_M:M\to C(M)$
is always a quasi-isomorphism, we can conclude that $M$ is closed if
and only if $\rho_M:M\to C(M)$ is an homotopy equivalence.
Using this characterization, by functoriality of $C(-)$, 
assuming that $\rho_X:X\to C(X)$ is an homotopy
equivalence it is clear that $\rho_X\ot id_V=\rho_{X\ot V}:X\ot V\to
C(X\ot V)=C(X)\ot V$ is also an homotopy equivalence, and consequently
$X\ot V$ is closed.

\bigskip

\noindent {\bf Notation:} We will say that a subcategory {\frak C}
of $\D^+(C)$ is closed by extensions of objects of $\Chain^+(k)$
if, whenever $X$ belongs to {\frak C} and $V\in\Chain^+(k)$,
then $X\ot V$ also belongs to {\frak C}.

\bigskip

As a particular case, taking $V=k^{(I)}$ with $I$ an arbitrary set,
the above corollary implies that if $X$ is closed,
then $X^{(I)}$ is closed as well. We can also prove the following
Lemma:

\begin{lem}\label{lemasd}
Let $C$ be a positively graded differential coalgebra and
$X\in\Chain^+(C)$ a direct summand of $Y\in\Chain^+(C)$
with $Y$ closed. Then $X$ is closed.
\end{lem}
\dem Let us call $Z$ a complement of $X$ in $Y$, i.e.
$X\oplus Z\cong Y$. Then we have a short exact sequence
$$0\to X\to \left((X\oplus Z)\oplus (X\oplus Z)\oplus\dots \right)\to
\left( Z \oplus (X\oplus Z)\oplus\dots \right) \to 0$$
The middle and right term (after rearranging parenthesis)
of this short sequence are isomorphic
to $Y^{(\nn)}$, then they are closed objects. Also the sequence
splits in $\Chain(C)$ (and then it obviously splits as a sequence
of graded $C$-comodules), then using Lemma \ref{lemaruso} the sequence
fits into a triangle, which proves that $X$ is closed.

\section{Derived equivalences}

\subsection[Equivalences induced by cotensor product]
{Equivalences induced by cotensor product and
inverse limits}

The main result of this section is the following Theorem,
that can be considered a partial dualization of
a result of B. Keller \cite{KE} 
on differential graded algebras.

\begin{teo}\label{cokeller}
Let $C$ and $D$ be two positively graded differential coalgebras,
$X\in\Chain(D\ot C^{op})$ with $C$ and $X$ bounded
complexes. Call $F:=X\Box_C^R-:\D^+(C)\to \D^+(D)$
and assume that $X$ is a closed object in $\Chain(D)$. Then
the following are equivalent:
\begin{enumerate}
\item $F$ is an equivalence.
\item \begin{enumerate}
	\item For all $n\in\ZZ$, $F$ induces an isomorphism
$\Hom_{\D(C)}(C,C[n])\cong \Hom_{\D(D)}({}_DX,{}_DX[n])$
	\item $F$ commutes with arbitrary products.
	\item The smallest triangulated full subcategory of $\D^+(D)$
which is stable under arbitrary products, extensions by 
objects of $\Chain^+(k)$ and contains ${}_DX$, contains ${}_DD$.
      \end{enumerate}
\end{enumerate}
\end{teo}

\dem $1\To 2)$. It is clear that $(a)$ and
$(b)$ are necessary conditions. If we consider the smallest
triangulated subcategory of $\D^+(C)$
which is stable under arbitrary products, extensions by 
objects of $\Chain^+(k)$ and contains ${}_CC$, then it
contains every standard resolution, so it is equivalent
to $\D^+(C)$. We recall that the standard resolution $C(M)$
is an inverse limit
in $\Chain(C)$ that can be viewed as a part of a triangle, the other
two objects in the triangle are products of $C_n(M)$. Each $C_n(M)$
can be constructed inductively with triangles by
``adding" objects like $C\ot(C^{\ot n-1}\ot M)$, and all this 
constructions involve products, triangles, and extensions of
$C$ by objects of $\Chain^+(k)$.

Since $F$ commutes with all the operations
mentioned above, and $F(C)={}_DX$, it follows that
the smallest
triangulated subcategory of $\D^+(D)$
which is stable under arbitrary products, extensions by 
objects of $\Chain^+(k)$ and contains ${}_DX$ is
equivalent
to $\D^+(D)$, in particular it contains $D$.

\medskip

$2\To 1).$ By the same arguments of the above paragraph,
$(b)+(c)$ implies that $F$ is quasi-surjective, because
$(c)$ states that ${}_DD\cong F(Y)$ for some $Y\in \D^+(C)$.
Then, given $N\in \D^+(D)$, its standard resolution
$D(N)=Tot(\oplus_{n\geq 1} D^{\ot n}\ot N,b',d)$ is isomorphic
to some object in the image of $F$, since the image of $F$ is stable
under products, triangles, extension by objects of $\Chain^+(k)$
and contains $D$.

The hardest part is to see that $F$ is fully faithful. We will 
need the notion of a certain kind of inverse limits, and the 
corresponding behaviour of these limits with respect to the
$\Hom$ functor. The proof of $2\To 1)$ will finish after Lemma
\ref{764} and \ref{765}.

\bigskip

Let $M\in\Chain^+(C)$ and call $M_{\leq n}$ the truncated
complex at degree $n$, we have that
$M=\stackunder{\lo n}{\lim}M_{\leq n}$ in $\Chain(k)$, 
and there is,
for all $X\in\Chain(k)$,
a natural injective morphism
$$\Hom_{Ch(k)}(M,X)\hookrightarrow
\ \stackunder{\to n}{\lim}\Hom_{Ch(k)}(M_{\leq n},X)$$
In case that $X_p=0$ for $p\gg 0$, the above morphism is clearly
an isomorphism.

The problem is that in general, $M_{\leq n}$ is not an 
object of $\Chain(C)$ (unless for example when $C$ is
a concentrated coalgebra). So we will define a class of inverse 
limits satisfying a kind of Mittag-Leffler condition, with 
analogous properties with respect to the $\Hom$ funtor, but
in the category $\Chain(C)$:

\begin{defi}
Let $\{U^{n+1}\to U^{n}\}_{n\in\nn}$ be an inverse system 
in $\Chain(C)$. We call the system {\bf locally finite} if,
given $j\in\NN$, there exists $n_0(j)$ such that the map
$U^{n+1}_{\leq j}\to U^{n}_{\leq j}$ is the identity
for all $n\geq n_0(j)$.
\end{defi}

\bigskip

\ex Let $C$ be a positively graded differential coalgebra
and $M\in\Chain^+(C)$, then $\{C_n(M)\}_{n\geq 1}=
\{(\oplus_{j=1}^nC^{\ot j}\ot M,b'+d)\}_{n\geq 1}$ 
is a locally finite inverse system.

\begin{lem}
Let $\{U^n\}_{n\in\nn}$ be a locally finite inverse system
in $\Chain^+(C)$ (i.e.
$U^n\in\Chain^+(C)$ and $U:=\stackunder{\lo n}{\lim}U^n\in \Chain^+(C)$)
and let $W\in\Chain^b(C)$. Then
there is a canonical isomorphism
$$\Hom_{Ch(C)}(U,W)\cong \ \stackunder{\to n}{\lim}\Hom_{Ch(C)}(U^n,W)$$
\end{lem}
\dem Let $P\in\NN$ be a number such that $W_p=0\ \forall \ \ |p|\geq P$.
Since $\{U^n\}_{n\geq 1}$ is a locally finite inverse system,
there exists $m_0(P)$ such that the maps
$(U^{n+1})_{\leq m}\to (U^{n})_{\leq m}$ are the identity for all
$m\geq m_0(P)$. We remark that if $f:U^n\to W$ is a map
in $\Chain(C)$ and $f_p$ denotes the $p$-th component of $f$,
then $f_p=0\ \forall \ \ |p|>P$. This remark, together with the 
locally finiteness condition implies that if
$\{f_n:U^n\to W\}$ is a system of maps, compatible with
the maps of the system $\{U^n\}$, then the map on the limit
is determined by $f^m:U^m\to W$ for any $m\geq m_0$,
then
$$\Hom_{Ch(C)}(U,W)=\Hom_{Ch(C)}(U^{m_0},W)=
\stackunder{\to n}{\lim}(\Hom_{Ch(C)}(U^n,W))$$

\begin{lem}\label{764}
Let $X\in\Chain^b(C)$, $F:\Chain^+(C)\to \Chain^+(D)$ a functor
commuting with locally finite inverse limits and
assume $F(X)\in\Chain^b(D)$. Let $V=\stackunder{\lo i}{\lim}V^i
\in\Chain^+(C)$ a limit of a locally inverse system in $\Chain^+(C)$
such that, for every $i\in\NN$, $F$ induces an isomorphism
$$F:\Hom_{\H(C)}(V^i,X)\to \Hom_{\H(C)}(F(V^i),F(X))$$
Then $F$ induces an isomorphism
$$F:\Hom_{\H(C)}(V,X)\to \Hom_{\H(C)}(F(V),F(X))$$
\end{lem} 
\dem Given two complexes $Y,\ Z\ \in \Chain(C)$, we denote
by $\H om_C(Y,Z)_*$ the $\H om$ complex (for details see for example
\cite{har}). This complex has in each degree
$$\H om_C(Y,Z)_n=\Hom_{C}^{\zz}(Y,Z[n])$$
i.e. graded $C$-colinear morphism from $Y$ to $Z[n]$.
The cohomology  of this complex in degree zero are $\ZZ$-graded
$C$-colinear maps which commute with differentials
(i.e. morphisms of $\Chain(C)$), modulo the homotopy equivalence
relation. We then obtain:
$$\begin{array}{rl}
\Hom_{\H(C)}(V,X)=&H^0(\H om_C(V,X))\\
=&H^0(\H om_C(\stackunder{\lo i}{\lim}V^i,X))\\
=&H^0(\stackunder{\to i}{\lim}\H om_C(V^i,X))\\
=&\stackunder{\to i}{\lim}H^0(\H om_C(V^i,X))\\
=&\stackunder{\to i}{\lim}\Hom_{\H(C)}(V^i,X)\\
=&\stackunder{\to i}{\lim}\Hom_{\H(D)}(F(V^i),F(X))\\
=&\stackunder{\to i}{\lim}H^0(\H om_{D}(F(V^i),F(X)))\\
=&H^0(\stackunder{\to i}{\lim}\H om_{D}(F(V^i),F(X)))\\
=&H^0(\H om_{D}(\stackunder{\lo i}{\lim}F(V^i),F(X)))\\
=&H^0(\H om_{D}(F(\stackunder{\lo i}{\lim}V^i),F(X)))\\
=&H^0(\H om_{D}(F(V),F(X)))\\
=&\Hom_{\H(D)}(F(V),F(X))
\end{array}$$

\bigskip

\begin{lem}\label{765}
Let $F:\Chain(C)\to\Chain(D)$ be a functor
and $\{U^n\}_{n\in\nn}$ an inverse system
in $\Chain(C)$ such that the sequences
$$0\to U\to \prod_n U^n \to \prod_n U^n\to 0$$
$$0\to F(U)\to \prod_n F(U^n) \to \prod_n F(U^n)\to 0$$
are triangles in $\H(C)$ (for example if the maps
$U^{n+1}\to U^n$ and
$F(U^{n+1})\to F(U^n)$ 
are split surjections
as morphisms of graded comodules)
where $U=\stackunder{\lo}{\lim}U^n$.
Assume that for every $n\in\NN$,
the functor $F$ induces an isomorphism
$F:\Hom_{\H(C)}(Z,U^n)\cong \Hom_{\H(D)}(F(Z),F(U^n))$,
then
$$F:\Hom_{\H(C)}(Z,U)\cong \Hom_{\H(D)}(F(Z),F(U))$$
\end{lem}

\dem we apply the $\Hom_{\H}$ functor to the triangles of the 
hypothesis and we obtain a morphism of long exact sequences:
{\footnotesize $$\diagram
\dots
\rto &\Hom_{\H(C)}(Z,U)\ddouble\rto &
\Hom_{\H(C)}(Z,\stackunder{n}{\prod} U^n) \rto \ddouble&
\Hom_{\H(C)}(Z,\stackunder{n}{\prod} U^n)\rto \ddouble& \dots\\
\dots
\rto &\Hom_{\H(C)}(Z,U)\rto \dto^{F}&
\stackunder{n}{\prod} \Hom_{\H(C)}(Z,U^n) \rto \ddouble&
\stackunder{n}{\prod} \Hom_{\H(C)}(Z, U^n)\rto \ddouble& \dots\\
\dots
\rto &\Hom_{\H(D)}(F(Z),F(U))\rto &
\stackunder{n}{\prod} \Hom_{\H(D)}(F(Z),F(U^n)) \rto &
\stackunder{n}{\prod} \Hom_{\H(D)}(F(Z),F(U^n))\rto & \dots\\
\enddiagram$$}
Then by the five Lemma,
$F:\Hom_{\H(C)}(Z,U)\cong \Hom_{\H(D)}(F(Z),F(U))$.

\bigskip

Now we study the behaviour of the locally finite inverse 
systems with respect
to the tensor and cotensor product:

\begin{lem}
Let $X\in\Chain^+(k)$ and $\{U^n\}_{n\geq 1}$ be a locally finite
inverse system in $\Chain^+(k)$, then
$$\stackunder{\lo}{\lim}(X\ot U^n)=X\ot U$$
\end{lem}
\dem We first remark that if $X$ and $Y$ are objects
of $\Chain^+(k)$, then only
a finite number of summands appear in
 each component of the tensor
product $(X\ot Y)_r=\oplus_{p+q=r}X_p\ot Y_q$, then
$X\ot Y=\stackunder{\lo j}{\lim}(X\ot Y_{\leq j})$. Now following
the definition of locally finite inverse system we have that

$$X\ot U=
\stackunder{\lo j}{\lim}(X\ot U_{\leq j})=
\stackunder{\lo j}{\lim}(X\ot \stackunder{\lo n}{\lim}U^n_{\leq j})=$$
$$=\stackunder{\lo j}{\lim}(X\ot U^{n_0(j)}_{\leq j})=
\stackunder{\lo j}{\lim}(\stackunder{\lo n}{\lim}  (X\ot U^n_{\leq j})=$$
$$=\stackunder{\lo n}{\lim}(\stackunder{\lo j}{\lim}X\ot U^n_{\leq j})=
\stackunder{\lo n}{\lim}( X\ot U^n)$$

\begin{coro}
Let ${}_DX_C\in\Chain^+(C\ot D^{op})$  and $\{U^n\}_{n\geq 1}$
be a locally finite inverse system in $\Chain^+(C)$, where $C$ and $D$ 
are positively graded differential coalgebras. Then 
$$\stackunder{\lo n}{\lim}(X\Box_CU^n)=
X\Box_C\stackunder{\lo n}{\lim}U^n$$
\end{coro}
\dem once we know that $\ot$ commutes with this particular class
of inverse limits, we only notice that $X\Box_C-$ is defined
as a kernel of the form
$$0\to X\Box_C U\to X\ot U\to X\ot C\ot U$$
and kernels commutes with arbitrary inverse limits.

\bigskip

Now we come back to the proof of Theorem \ref{cokeller}, we have to
prove that if the functor ${}_DX_C\Box^R-$ verifies $(a)$ 
$(b)$ and $(c)$ then it is fully
faithful.

Let us call ${\frak C}_1$ the full subcategory of
$\D^+(C)$ whose objects are complexes $U$ such that,
for all $n\in\ZZ$,
$F$ induces an isomorphism
$$F:\Hom_{\D(C)}(U,C[n])\to \Hom_{\D(D)}(F(U),{}_DX[n])$$
It is clear that ${\frak C}_1$ is a triangulated subcategory
of $\D^+(C)$ and by $(a)$ it contains ${}_CC$. Also by 
Lemma \ref{764} it is stable under locally finite
inverse limits, then ${\frak C}_1$ coincides with $\D^+(C)$.

Let us now consider ${\frak C}_2$ the full subcategory
$\D^+(C)$ whose objects are complexes $V$ such that,
for all $U\in\D^+(C)$,
$F$ induces an isomorphism
$$F:\Hom_{\D(C)}(U,V)\to \Hom_{\D(D)}(F(U),F(V))$$
It is clear that ${\frak C}_2$ is a triangulated subcategory
of $\D^+(C)$, stable by direct summands, and by $(b)$ 
${\frak C}_2$ is stable under arbitrary products. By the above discussion
${\frak C}_2$ contains $C$, in order to see that
${\frak C}_2$ is equivalent to $\D^+(C)$ it is enough to
see that $C\ot V \in {\frak C}_2$ for all $V\in \Chain^+(k)$.

Since every object $V\in\Chain^+(k)$ is a locally finite inverse limit
of objects of
$\Chain^b(k)$, by Lemma \ref{765} one can assume that $V$ is bounded.
Assume that $V$ is of the form
$$\diagram 
\dots\rto&
0\rto&
V_{n}\rto^{d_n}&
V_{n+1}\rto^{d_{n+1}}&
\dots\rto^{d_{m-1}}&
V_{m}\rto&0\rto&\dots\enddiagram$$
Let us call $V_{\leq m-1}$ the complex which has the same components
as $V$ on every degree except in degree $m$, and the same differentials
(except $d_{m-1}$), then $d_{m-1}:V_{\leq m-1}\to V_{m}$ is
a morphism of complexes, and $V\cong Co(d_{m-1})$. Now the funtor
$C\ot-:\D^+(k)\to \D^+(C)$ commutes with mapping cones, then inductively,
since ${\frak C}_2$ is stable under mapping cones, 
it is enough to see that $C^{(I)}$ belongs to
${\frak C}_2$, for every set $I$, and this is easy, because
$k^{(I)}$ is (as $k$-vector space) a direct summand
of $k^I$, so $C^{(I)}$ is a direct summand of
$C\ot k^I=C^I$ (the product in the category $\Chain^+(C)$),
then $C^{(I)}$ belongs to ${\frak C}_2$ because
${\frak C}_2$ is stable under products and direct summands.

\subsection{On derived equivalences of concentrated coalgebras}

Let us consider the case when $C$ is a usual coalgebra and we 
view it as a differential graded one with trivial differential
graded structure. 

Free comodules are closed objects because they are of the form
$C^{(I)}=C\ot k^{(I)}$, and direct summands of closed objects
are also closed, then injective comodules are closed.
Being closed is also stable by triangles, then inductively 
one can easily prove that a bounded complex with injective 
components is closed (for example the complex 
$\dots 0\to X_n\to^d X_{n+1}\to 0\to \dots$
is the mapping cone of the map $d:(\dots 0\to X_n\to 0\to \dots)
\to (\dots 0\to X_{n+1}\to 0\to \dots)$). Also if
$M\in\Chain^+(C)$ is a complex with injective components, 
then $M=\stackunder{\lo n}{\lim}M_{\leq n}$
where $M_{\leq n}$ is the complex truncated at degree $n$, and
the surjections $M_{\leq n+1}\to M_{\leq n}$ clearly split
as morphisms of graded $C$-comodules, then by Lemma \ref{tri}
$M$ is a closed object.

On the other hand, if $M$ is any object in $\Chain^+(C)$, the 
standard resolution $C(M)$ has, in each degree, free 
$C$-comodules. We characterize then, in this special case, 
the class of closed objects as those ones homotopy equivalent to
(left bounded) complexes with injective components.

We can state the following proposition:

\begin{prop}\label{propb}
Let $C$ and $D$ be two concentrated coalgebras. 
Let $\H_I$ (resp. $\H_I^+$, $\H_I^b$) be the subcategory of 
$\H$ (resp. $\H^+$, $\H^b$) consisting of complexes with
injective components. 
Then $\D^+(C)\cong \D^+(D)$ (as triangulated categories)
if and only if  $\H_I^+(C)\cong \H_I^+(D)$, and any 
triangulated equivalence 
$\H_I^+(C)\cong \H_I^+(D)$ restricts to an equivalence
$\H_I^b(C)\cong \H_I^b(D)$.
\end{prop}

\dem we know in general that for any positively 
graded differential coalgebra, the category 
$\D^+$ is equivalent to
the subcategory of $\H^+$ consisting of closed objects,
denoted by $\H^+_c$. The discussion
above proves that when the coalgebra is concentrated,
the inclusion $\H_I^+\to \H^+_c$ is an equivalence, then 
for $C$ and $D$ as in the hypothesis we have that
$\D^+(C)\cong \H_I^+(C)$ and $\D^+(D)\cong \H_I^+(D)$, 
and the first assertion is clear.

\medskip

In order to see that any triangulated equivalence
$\H_I^+(C)\cong\H_I^+(D)$ restricts to an equivalence
$\H_I^b(C)\cong\H_I^b(D)$ we will need the following characterization
of bounded complexes:

\begin{lem}
Let $C$ be a concentrated coalgebra and $X\in\H^+(C)$. Then
$X\in\H^b(C)$ if and only if, for all $Y\in\H^+(C)$ there
exists $n_0\in\NN$ such that $\Hom_{\H(C)}(X[n],Y)=0$
for all $n\leq n_0$.
\end{lem}
\dem Let us first see that the condition is necessary. 

Assume $X\in\H^b(C)$, then $\exists \ m\in\NN$ such that
$X_p=0\ \ \forall \ p\ /\ |p| >m$. Now given $Y\in\Chain^+(C)$
there exists $m'\in\ZZ$ such that $Y_p=0\ \forall \ p<m'$.
If we take $n\leq m-m'$ then $\Hom_{Ch(C)}(X[n],Y)=0$ which
obviously implies $\Hom_{\H(C)}(X[n],Y)=0$.

\medskip

Let us now see that the condition is sufficient:

Let $X\in\H^+(C)$ and consider $Y=C$. By hypothesis there exists
$n_0\in\NN$ such that $\Hom_{\H(C)}(X[n],C)=0$ for all $n\leq n_0$.
By Proposition \ref{adjuncion},
$$\Hom_{\H(C)}(X[n],C)=
\Hom_{\H(k)}(X[n],k)=X_{-n}^*$$
then $X_{n}=0$ for all $n\geq n_0$, i.e. $X\in \H^-(C)$, then
$X\in\H^b(C)$.

\bigskip

\rems\begin{enumerate}
\item After this Lemma, the proof of Proposition \ref{propb}
is completed.
\item  Proposition \ref{propb}
tells us that, given an equivalence
$F:\D^+(C)\to \D^+(D)$ where $C$ and $D$ are concentrated coalgebras,
one can always assume that $F(C)$ is quasi-isomorphic to a bounded
complex with $D$-injective components, so the assumption
on the $X\in\Chain(D\ot C^{op})$ 
(about being bounded and $D$-closed)
in Theorem \ref{cokeller} is superfluous when $C$ and
$D$ are concentrated
coalgebras.
\item Let $F:\H_I^+(C)\to\H_I^+(D)$ a functor commuting with 
products and extensions by objects of $\Chain^+(k)$,
where $C$ and $D$ are concentrated coalgebras. Suppose
that $F$ restricts to a functor $\H_I^b(C)\to\H_I^b(D)$, and
assume that this restriction is an equivalence, then
re-writing the proof of Theorem \ref{cokeller}
we have that the original funtor $F$ is an equivalence. In
that sense, derived equivalences between concentrated
coalgebras can be ``checked" looking at the category $\H_I^b$.
\end{enumerate}

\section{Cotilting theory}

If we specialize Theorem \ref{cokeller} to the 
non differential graded case when
$C$ and $D$ are usual coalgebras, and $X$ is a $D$-$C$-bicomodule,
condition $(a)$ is a vanishing condition of the $Ext$ groups, and
condition $(b)$ is always satisfied when $X$ is quasi-finite
(because in that case $X\Box_C-$ admits a left adjoint).
Condition $(c)$ of \ref{cokeller} can be verified if
for example $D$ fits into some exact sequence where
the other components of the sequence are obtained by
some operations on $T$.

\begin{defi}\label{deftilting}
Let $C$ and $D$ be two (concentrated) coalgebras and ${}_DT_C$
a $D$-$C$-bicomodule. We will call $T$ a {\bf cotilting comodule} if


\begin{enumerate}
\item ${}_DT$ is quasi-finite, $e_D(T)\cong C$ and $e_C(T)\cong D$.
\item $Ext^n_D(T,T)=0\ \forall \ n\geq 1$.
\item There exists an exact sequence of type
$$0\to {}_DT_n \to \dots\to {}_DT_0 \to {}_DD\to 0$$
with $T_i\in Add(T)$.
\item There exists an exact sequence of $D$-$C$-bicomodules
$$0\to {}_DT_C \to I_0\to \dots\to I_r \to 0$$
with $I_i$ injective and quasifinite as $D$-comodules.
\end{enumerate}
\end{defi}

The main result of this section is the following:

\begin{teo}\label{teocotilting}
Let $C$ and $D$ be two (concentrated) coalgebras admiting a cotilting bicomodule
${}_DT_C$, and call $Q:=h_D(T,D)$. Then there are isomorphisms 
in $\D(D^e)$ and $\D(C^e)$:
$$T\Box_C^R Q\cong D$$
$$Q\Box_D^R T\cong C$$
In particular, $\D^+(C)\cong \D^+(D)$, as triangulated categories.
\end{teo}

Since we will not only prove that $T\Box_C^R-$ is an 
equivalence, but that its quasi-inverse is also a derived 
cotensor product, the proof  of the above Theorem is a little 
more complicated than just applying Theorem 
\ref{cokeller}. Nevertheless,
we begin with a Lemma in the direction of  condition $(c)$ of
\ref{cokeller}.

\begin{lem}\label{3.3}
Let $C$ and $D$ be two concentrated coalgebras, ${}_DT_C$ a bicomodule
which is quasi-finite as $D$-comodule, then  for all $X\in\Chain(C)$ and
$Y\in\Chain(D)$ we have natural isomorphisms:
$$\Hom_{Ch(C)}(h_D(T,Y),X)\cong\Hom_{Ch(D)}(Y,T\Box_CX)$$
$$\Hom_{\H(C)}(h_D(T,Y),X)\cong\Hom_{\H(D)}(Y,T\Box_CX)$$
where the differential of $h_D(T,Y)$
and  $T\Box_CX$ are respectively  $h_D(T,d_Y)$ and
$id\Box_Cd_X$.
\end{lem}
\dem we will make use of the complex $\H om$. We recall that the 
component in degree $p$ of $\H om$ is the set of homogeneous 
$C$-colinear maps of degree $p$,
more precisely, we have that
$$\H om_C(h_D(T,Y),X)_p=\prod_{n\in\zz}\Hom_C(h_D(T,Y_n),X_n[p])\cong$$
$$\cong \prod_{n\in\zz}\Hom_D(Y_n,T\Box_CX_n[p])=\H om_D(Y,T\Box_CX)_p$$
The fact that this isomorphism commutes with the differential
in $\H om$ comes from the naturality of the isomorphism applied 
to $d_X^n:X_n\to X_{n+1}$
and $d_Y^n:Y_n\to Y_{n+1}$. 

The complexes $\H om$ being isomorphic, their cocycles 
in degree zero are isomorphic and they have 
the same cohomology in degree zero, and the proof 
of the Lemma is completed.

\begin{coro}\label{coroprod}
Under the same hypothesis of Lemma \ref{3.3}, 
$T\Box_C^R-:\D^+(C)\to \D^+(D)$ commutes
with products.
\end{coro}
\dem the functor $T\Box_C-:\H(C)\to \H(D)$ commutes with products 
because it admits a left adjoint. Let $\{M_i\}_{i\in I}$ a family of objects 
in $\Chain^+(C)$ such that $\prod_{i\in I}M_i\in\Chain^+(C)$, in order to
compute $T\Box_C^R\left(\prod_{i\in I}M_i\right)$ we need to find a 
closed object in $\Chain^+(C)$ quasi-isomorphic to $\prod_{i\in I}M_i$.
One option is to consider
$C\left(\prod_{i\in I}M_i\right)$, but we consider instead
$\prod_{i\in I}C(M_i)$, which is closed because 
it is a product of closed objects, and
it is also  quasi-isomorphic to
$\prod_{i\in I}M_i$ because each $M_i$ is $k$-homotopically 
equivalent to $C(M_i)$, then the product
of the $M_i$'s is $k$-homotopically equivalent 
to the product of the $C(M_i)$'s. We obtain then
$$T\Box_C^R\left( \prod_{i\in I}M_i \right)=
T\Box_C\left( \prod_{i\in I}C(M_i) \right)=$$
$$= \prod_{i\in I}T\Box_CC(M_i)
= \prod_{i\in I}T\Box_C^RM_i$$
\medskip

We can now prove a different version of theorem \ref{teocotilting}:

\begin{teo}\label{Tbox}
Let $C$ and $D$ be two concentrated coalgebras and ${}_DT_C$ a 
bicomodule satisfying conditions 1. 2. and 3. of  the definition 
of a cotilting bicomodule, then $T\Box_C^R-:\D^+(C)\to \D^+(D)$ is
an equivalence of triangulated categories.
\end{teo}
\dem we use the characterization given in Theorem \ref{cokeller}.
\begin{itemize}
\item[(a)] $\Hom_{\D(C)}(C,C[n])=Ext^n_C(C,C)$ which is zero 
unless $n=0$, in that case  $Ext^0_C(C,C) \cong 
Com_C(C,C)=C^*$. On the other
hand, $\Hom_{\D(D)}(T,T[n])=Ext^n_C(T,T)$ which by hypothesis
is zero for $n\geq 1$, and $Ext^0_C(T,T) = Com_D(T,T) = (e_D(T))^*
\cong C^*$.
\item[(b)] it is Corollary \ref{coroprod}.
\item[(c)] if a triangulated subcategory of $\D^+(D)$ contains
$T$ and is stable by extensions by objects of $\Chain^+(k)$, then
it contains $T^{(I)}$ and (with same proof of Lemma \ref{lemasd})
it contains its direct summands, then it contains $Add(T)$. Now consider
the exact sequence
$$0\to {}_DT_n \to \dots\to {}_DT_0 \to {}_DD\to 0$$
with $T_i\in Add(T)$. Inductively, the complex
$$\diagram 
0\rto &T_n\rto^{d_n}&\dots \rto^{d_i}&T_{i-1}\rto & 0
\enddiagram$$ is the mapping cone of
the map $d_i$:
$$\diagram 
0\rto &T_n\rto^{d_n}&\dots \rto    &
		   T_{i+1}\rto^{d_{i+1}}&T_{i}\dto^{d_i}\rto & 0\\
0\rto &0   \rto            &\dots \rto    &
 0           \rto                   &T_{i-1}\rto    &0 \\
\enddiagram$$
viewed as a map of complexes.
So the complex 
$0\to {}_DT_n \to \dots\to {}_DT_0 \to 0$ belongs to the triangulated
subcategory of $\D^+(D)$ stable under extensions by objects
of  $\Chain^+(k)$, and it is quasi-isomorphic to $D$.
\end{itemize}

Now coming back to Theorem \ref{teocotilting}, we will prove
a Lemma of comodule theory which may be probably considered part
of the folklore of the theory, we include the proof for completeness:

\begin{lem}
Let $C$ and $D$ be two (concentrated) coalgebras and
${}_DI_C$ a bicomodule 
such that, viewed as $D$-comodule, is injective
and quasi-finite. Then
for all $D$-comodule $M$ there is an isomorphism of $C$-comodules
$$h_D(I,M)\cong h_D(I,D)\Box_D M$$
\end{lem}
\dem The functor $h_D(I,-)$ is left adjoint to the functor
$I\Box_C-$, then it commutes with direct sums, in particular, 
if $W$ is any $k$-vector space, we have that
$h_D(I,X\ot W)=h_D(I,X)\ot W$ for every $D$-comodule $X$.

On the other hand, the functor $h_D(I,-)^*=Com_D(-,I)$, 
and since $I$ is assumed to be an injective $D$-comodule, then
 $h_D(I,-)$ is exact. Now let $M$ be any $D$-comodule and apply the
funtor $h_D(I,-)$ to the exact sequence
$$
\diagram 
0\rto & M \rto^{\rho_M}& D\ot M\rto^{b'}&D\ot D\ot M
\enddiagram$$
By the above discussion we have the identifications
$$\diagram
0\rto & h_D(I,M) \rto &h_D(I,D\ot M)\ddouble \rto &
 h_D(I,D\ot D\ot M)\ddouble\\
0\rto & h_D(I,D)\Box_D M \rto &h_D(I,D)\ot M \rto & 
h_D(I,D)\ot D\ot M\\
\enddiagram$$
giving $ h_D(I,M)= h_D(I,D)\Box_D M$.


\begin{lem}
Let $C$ and $D$ be two positively
graded differential coalgebras and
$T\in\Chain(D\ot C^{op})$ such that, for every $Y\in\Chain(C)$
there exists an object $h_D(T,Y)$ in $\Chain(D)$ which is functorial
in $Y$, and an isomorphism of complexes
$\H om_D(h_D(T,Y),X)\cong \H om_C(Y,T\Box_CX)$ which is natural in $X$.
Then $e_D(T):=h_D(T,T)$ is a differential graded coalgebra.
\end{lem}
\dem We proceed in the same way as in the non-differential graded case:
from the identity map of $e_D(T)$, which is an element
of $\H om_D(e_D(T),e_D(T))_0$ one gets an element
of $\H om_C(T,T\Box_Ce_D(T))_0\subset
\H om_C(T,T\ot e_D(T))_0$, and iterating the procedure
one get an homogeneous map from $T$ into
$T\Box_Ce_D(T)\ot e_D(T)$. Now considering $X=
e_D(T)\ot e_D(T)$ and $Y=T$, we get an element
in
$\H om_D(e_D(T),e_D(T)\ot e_D(T))_0$.

The coassociativity of $e_D(T)$ is proved exactly in the
same way as in the non-differential graded case. We will
only remark that 
the colagebra structure and the differential graded structure
fit together because
the graded
dual object is an associative differential algebra, namely
$$e_D(T)^*=\H om_k(h_D(T,T),k)\cong$$
$$\cong \H om_C(h_D(T,T),C)\cong
\H om_D(T,T\Box_C C)\cong \H om_D(T,T)$$

\bigskip

\rems 1. the hypothesis of this Lemma include the case when
$T$ is a bicomodule over two concentrated coalgebras $C$ and $D$
and $T$ is quasi-finite viewed as $D$-comodule.

\medskip

\noindent 2. If $T$ and $T'$ are complexes 
satisfying the hypothesis of this Lemma and $f:T\to T'$
is a morphism in $\Chain(C)$, then
$Co(f)$ also satisfies the same hypothesis. 

\medskip

By the second remark, the class of objects satisfying an adjoint-type
hypothesis like the above one is closed under finite direct sums and 
shifting of degree. This implies that if $C$ and $D$ 
are two concentrated coalgebras, then any bounded complex
of $D$-$C$-bicomodules with $D$-quasi-finite components
satisfies the adjoint-type hypothesis.

\bigskip

Consider now $C$ and $D$ two concentrated
coalgebras admiting a cotilting bicomodule ${}_DT_C$.
By condition 4. of Definition \ref{deftilting}, there
exists an exact sequence of $D$-$C$-bicomodules
$$0\to {}_DT_C \to I_0\to \dots\to I_r \to 0$$
with $I_i$ injective and quasi-finite as $D$-comodules. 
\begin{lem}
Keeping notations, 
let us
call $T_*$ the complex
$0 \to I_0\to \dots\to I_r \to 0$. Then
the differential graded coalgebra $h_D(T_*,T_*)$
is quasi-isomorphic to $C$ as $C$-bicomodule.
\end{lem}
\dem it is enough to see that the graded dual 
algebras are quasi-isomorphic:

$$H^n\left(
\H om_D(T_*,T_*)
\right)=
H^0\left(
\H om_D(T_*,T_*[n])
\right)=
\Hom_{\H(D)}(T_*,T_*[n])$$
But $T_*$ is a bounded complex of $D$-injective objects, then
it is a closed object and so
$$\Hom_{\H(D)}(T_*,T_*[n])=
\Hom_{\D(D)}(T_*,T_*[n])\cong$$
$$\cong \Hom_{\D(D)}(T,T[n])\cong
Ext^n_D(T,T)$$
and $Ext^n_D(T,T)=0$ unless $n=0$, and in this case
$Ext^0_D(T,T)=Com_D(T,T)=e_D(T)^*\cong C^*$.

\bigskip

\begin{coro}\label{coro1}
Let $C$ and $D$ be two (concentrated) coalgebras and ${}_DT_C$ a
cotilting bicomodule. Let us call $T_*$ the complex
$0\to  I_0\to \dots\to I_r \to 0$ 
of the definition of cotilting comodule,
which is quasi-isomorphic to ${}_DT_C$.
Then there is an isomorphism
in $\D(C^e)$:
$$h_D(T_*,D)\Box_D T_*\cong C$$
\end{coro}
\dem Since $h_D(T_*,D)_n=h_D(I_{-n},D)$, we have that
$$(h_D(I_*,D)\Box_D T_*)_n=
\bigoplus_{p+q=n}h_D(I_{-p},D)\Box_D I_q=
\bigoplus_{p+q=n}h_D(I_{-p},I_q)$$
Then $h_D(T_*,D)\Box_D T_*=h_D(T_*,T_*)$ which, by the
two Lemmae above, is quasi-isomorphic to
$C$.

\bigskip

We now prove our last Lemma needed  for the proof ot
Theorem \ref{teocotilting}:

\begin{lem}\label{extension}
Let $C$ and $D$ be two differential positively graded coalgebras
and $F:\D^+(C)\to \D^+(D)$ an equivalence with inverse
$G:\D^+(D)\to \D^+(C)$. Then, for any positively graded differential 
coalgebra $E$, one can extend the functors
$F$ and $G$ in order to have an equivalence
between $\D^+(C\ot E)$ and $\D^+(D\ot E)$.
\end{lem}
\dem  Let $X\in \Chain^+(C\ot E)$, then using the standard resolution
of $X$ with respect to $E$, there is a $C\ot E$-colinear quasi-isomorphism
$${}_{C,E}X\to Tot\left(
\bigoplus_{n> 0}{}_EE^{\ot n}\ot {}_CX, b'_E,d
\right)$$
We define then $\wh{F}(X):=
Tot\left(
\bigoplus_{n> 0}{}_EE^{\ot n}\ot {}_DF(X), b'_E,d
\right)$. We remark that we can apply $F$ because the 
$E$-structure map of the above complex 
is a $C$-colinear map. In an analogous way we define
$\wh{G}$, let us see that they are inverse to each other:

$$\wh{G}(\wh{F}(X))=
Tot\left(\bigoplus_{n> 0}
{}_EE^{\ot n}\ot {}_CG(\wh{F}(X))
\right) = Tot
\left(\bigoplus_{n>0} {}_EE^{\ot n}\ot {}_CG
      \left(Tot
            \left(\oplus_{n> 0}E^{\ot n}\ot {}_DF(X)
                  \right)
            \right)
      \right)$$
But $Tot\left(\oplus_{n> 0}E^{\ot n}\ot {}_DF(X),b'_E,d\right)$ 
is quasi-isomorphic
to $F(X)$ as $D$-comodules, then
$G\left(
Tot\left(
\bigoplus_{n> 0}E^{\ot n}\ot {}_DF(X)
\right)\right)$ is isomorphic to $G(F(X))$, obtaining
$$\wh{G}(\wh{F}(X))\cong
Tot\left(\bigoplus_{n> 0}{}_EE^{\ot n}\ot {}_CG(F(X))\right)
\cong Tot\left(\bigoplus_{n> 0}{}_EE^{\ot n}\ot {}_CX\right)\cong
{}_{C,E}X$$
The other isomorphism is proved identically.

\bigskip

As a Corollary, we can write down a complete proof of
Theorem \ref{teocotilting}:

\medskip

After Corollary \ref{coro1} we know that 
$h_D(T_*,D)\Box_D T_*\cong C$. We have then isomorphisms
in $\D(C^e)$:
$$C\cong h_D(T_*,D)\Box_D T_*\cong 
h_D(T,D)\Box_D T_*\cong  h_D(T,D)\Box_D^R T$$
We also know that $T\Box_C^R-$ is an equivalence, let us
call $F:=T\Box_C^R-$ and consider $\wh{F}:\D^+(C\ot D^{op})\to
\D^+(D^e)$, the extension given in Lemma \ref{extension}. Since
$\wh{F}$ is an equivalence, there must exist an object
$S\in \Chain(C\ot D^{op})$ such that $\wh{F}(S)\cong {}_DD_D$.
We have then the following isomorphisms in $\D(D^e)$:

$${}_DD_D\cong\wh{F}(S)=
Tot\left(
\bigoplus_{n>0} D^{\ot n}_D \ot {}_DT\Box_C^RS,b'_D,d
\right)=$$
$$=Tot\left(
\bigoplus_{n>0} D^{\ot n}_D \ot {}_DT_*\Box_CS,b'_D,d
\right)\cong
{}_DT_*\Box_CS_D$$
Now we recall that $h_D(T,D)\Box_D^RT\cong{}_CC_C$, and this
implies that $S$ must be isomorphic (in $\D(D^e)$) to
$h_D(T,D)$, as it is shown below:

$$S\cong C\Box_C^RS\cong h_D(T,D)\Box_D^RT\Box_C^RS\cong$$
$$\cong h_D(T,D)\Box_D^RD\cong  h_D(T,D)$$

\rems 1. If $C$ and $D$ are two (concentrated) coalgebras which are
Morita - Takeuchi equivalent (i.e. the category of $C$-comodules
is equivalent to the category of $D$-comodules) then their derived
categories are equivalent, and we know after \cite{tak} that there
exists an injective cogenerator ${}_DI_C$ such that
$I\Box_Ch_D(I,D)\cong D$ and
$h_D(I,D)\Box_DI\cong C$ (isomorphisms of bicomodules). This 
implies that the notion of cotilting bicomodule generalizes
the Morita - Takeuchi equivalence relation, and the generalization
is strict, because it is enough to take a cotilting bicomodule
which is not injective and so it can never give an equivalence
at the level of comodule categories, nevertheless it gives
an equivalence at the level of derived categories.

\medskip 

\noindent 2. If $C$ and $D$ are two concentrated coalgebras
and $F:\D^+(C)\to \D^+(D)$ is an equivalence of triangulated 
categories, one can ask if there exists an object ${}_DX_C$
in $\Chain^+(D\ot C^{op})$ such that $X\Box_C^R-:\D^+(C)\to
\D^+(D)$ is an equivalence. If this were the case, it is
clear that ${}_DX$ is the value of
$X\Box_C^R-$ in $C$, if we put ${}_DX=F(C)$ then one can 
suppose that $X$ is a bounded
complex with $D$-injective 
components (Proposition \ref{propb}). But in order to
have the $C$-structure, one should need to equip $F(C)$
with a structure of $C^{op}$-comodule. In general,
if some $X\in \Chain^+(D)$ admits a coendomorphisms
coalgebra $\E:=h_{D}(X,X)$, then $X$ becomes a
$D\ot\E^{op}$-comodule in a canonical way. But the
condition of admiting a coendomorphism coalgebra is quite
restrictive, for example this is the case if $X$ is
$D$-quasi-finite, and we do not know if ``quasi-finiteness"
is a property preserved by derived equivalences. This 
discussion is the reason  that led us to include 
condition 4. in the definition of a cotilting comodule $T$.

\section{A quasi-isomorphism, an example of derived equivalence}
Let $C$ and $D$ be two positively graded differential
coalgebras and $f:C\to D$ be a quasi-isomorphism of
differential graded coalgebras. Composing $f$ with the comultiplication
of $C$
$$\diagram 
C\rto^{\Delta}&C\ot C\rto^{id\ot f}&C\ot D
\enddiagram$$
we can consider $C$ as an object in $\Chain(C\ot D^{op})$. We will
denote $C_{f}$ this $C$-$D$-bicomodule, and similarly
${}_fC$ and ${}_fC_f$. The main result of this section is that the derived
category associated to $C$ is equivalent to
the derived category associated to $D$, more precisely,
we have the following result:

\begin{prop}\label{propquis}
With the above notations, there are isomorphisms
in $\D(D^e)$ and $\D(C^e)$:
$${}_fC\Box_C^RC_f\cong D$$
$$C_{f}\Box_D^R{}_fC\cong C$$
\end{prop}
\dem The easiest part is to show that
${}_fC\Box_C^RC_f\cong D$. It is clear that $C_{f}$ is a closed object
in $\Chain(C)$, then
${}_fC\Box_C^RC_f=
{}_fC\Box_CC_f={}_fC_f$. But by hypothesis $f:C\to D$ is a 
quasi-isomorphism of coalgebras, then
$f:{}_fC_f\to D$ is a quasi-isomorphism of $D$-bicomodules, hence
an isomorphism in the derived category.

\medskip

In order to compute $C_f\Box_D^R{}_fC$, 
we have to find an object in
$\Chain(C\ot D^{op})$ quasi-isomorphic to
$C_f$ and closed as right $D$-comodule. We consider the standard
resolution of $C_f$ as right $D$-comodule, then

$$C_f\Box_D^R{}_fC=Tot(\oplus_{n\geq 1}C\ot D^{\ot n},d,b')
\Box_D{}_fC\cong$$
$$\cong Tot(\oplus_{n\geq 1}C\ot D^{\ot n-1}\ot C,d,b')$$
On the other hand, using the standard resolution of $C$ as $C$-comodule, we
have that $C$ is quasi-isomorphic (as $C$-bicomodule) to
$Tot(\oplus_{n\geq 1}C^{\ot n+1},d,b')$. 
A $C^e$-colinear morphism of complexes 
$$Tot(\oplus_{n\geq 1}C^{\ot n+1},d,b')\to
Tot(\oplus_{n\geq 1}C\ot D^{\ot n-1}\ot C,d,b')$$
is obtained applying $f$ in the middle terms.
Next we will see that the above map is
a quasi-isomorphism.

We know by the K\"unneth formula that 
$id\ot f^{\ot r}\ot id: (C^{\ot r+2},d)\to (C\ot D^{\ot r}\ot C,d)$
is a quasi-isomorphism, looking at the double complexes
$$\diagram
\vdots&\vdots&\vdots&\\
(C\ot C)_2\rto^{b'}\uto^d&(C\ot C\ot C)_{2}\rto^{b'}\uto^d&
(C\ot C^{\ot 2}\ot C)_{2}\rto^{\ \ b'}\uto^d&\dots\\
(C\ot C)_1\rto^{b'}\uto^d&(C\ot C\ot C)_{1}\rto^{b'}\uto^d&
(C\ot C^{\ot 2}\ot C)_{1}\rto^{\ \ b'}\uto^d&\dots\\
(C\ot C)_0\rto^{b'}\uto^d&(C\ot C\ot C)_{0}\rto^{b'}\uto^d&
(C\ot C^{\ot 2}\ot C)_{0}\rto^{\ \ b'}\uto^d&\dots\\
\enddiagram$$
and
$$\diagram
\vdots&\vdots&\vdots&\\
(C\ot C)_2\rto^{b'}\uto^d&(C\ot D\ot C)_{2}\rto^{b'}\uto^d&
(C\ot D^{\ot 2}\ot C)_{2}\rto^{\ \ b'}\uto^d&\dots\\
(C\ot C)_1\rto^{b'}\uto^d&(C\ot D\ot C)_{1}\rto^{b'}\uto^d&
(C\ot D^{\ot 2}\ot C)_{1}\rto^{\ \ b'}\uto^d&\dots\\
(C\ot C)_0\rto^{b'}\uto^d&(C\ot D\ot C)_{0}\rto^{b'}\uto^d&
(C\ot D^{\ot 2}\ot C)_{0}\rto^{\ \ b'}\uto^d&\dots\\
\enddiagram$$
we see that our map is a quasi-isomorphism on the columns,
so filtering by the columns, by a standard spectral sequence
argument we have that the total complexes are quasi-isomorphic.

%

\section[$\Hoch^*$, $H^*$ and $HC^*$ for D.G. coalgebras]
{$\Hoch^*$, $H^*$ and $HC^*$ for differential graded coalgebras}

Given a (concentrated) $k$-coalgebra $C$ and a bicomodule
$M$, Doi defined in \cite{doi} two cohomology theories
$\Hoch^*(M,C)$ and $H^*(M,C)$ which play the role of Hochschild
homology and cohomology, respectively. They are defined in terms
of standard complexes, but also have an interpretation as
derived functors, namely 
$$Hoch^*(M,C)=Cotor^*_{C^e}(M,C)=H^*(M\Box_{C^e}^RC)$$
$$H^*(M,C)=Ext^*_{C^e}(M,C)$$
In \cite{hc}, A. Solotar and myself defined a cyclic cohomology theory for 
coalgebras (denoted by $HC^*$) and proved some of its fundamental 
properties, 
including Morita - Takeuchi invariance (see \cite{fs} for the
invariance of $\Hoch^*$ and $H^*$ and \cite{hc} for the invariance
of $HC^*$ and a more general proof of the invariance of 
$\Hoch^*$). More precisely
we proved that
given two $k$-coalgebras $C$ and $D$ such that there exists a $k$-linear
equivalence $F:C$-comod$\to D$-comod, then there exists an equivalence
$\wh{F}:C^e$-comod$\to D^e$-comod with $\wh{F}(C)=D$
such that, given a  $C$-bicomodule $M$
$$\Hoch^*(M,C)\cong \Hoch^*(\wh{F}(M),D)$$
$$H^*(M,C)\cong H^*(\wh{F}(M),D)$$
$$HC^*(C)\cong HC^*(D)$$

The purpose of this section is to extend the definition of the
cohomology theories to the differential graded case, generalizing 
also the invariance results to derived equivalences. 

\bigskip

Let $C$ be a differential graded coalgebra and consider, for
each $n\in\NN_0$, the vector space $C^{\ot n +1}$. 
The following 
natural operators are defined on them:
\begin{itemize}
\item The differential $d:C^{\ot n+1}\to C^{\ot n+1}$
$$(c_0,\dots,c_n)\mapsto (-1)^n\left(
\sum_{k=0}^n(-1)^{|c_0|+\dots+|c_{k-1}|}(c_0,\dots,d_C(c_k),\dots,c_n)
\right)$$
\item The cyclic operator $T$
$$(c_0,\dots,c_n)\mapsto (-1)^n(-1)^{|c_0|.|(c_1,\dots,c_n)|}
(c_1,\dots,c_n,c_0)$$
where $|(c_1,\dots,c_n)|$ is the
standard degree on the tensor product:
$|(c_1,\dots,c_n)|=\sum_{i=1}^n|c_i|$.
\item $\Delta_i:C^{\ot n+1}\to C^{\ot n+2}$, $i=0,\dots, n$
$$(c_0,\dots,c_n)\mapsto (c_0,\dots,\Delta(c_i),\dots ,c_n)$$
and $\Delta_{n+1}:=(-1)^{n+1}T\Delta_0$.
\item The differential $b':=\sum_{i=0}^n(-1)^i\Delta_i$
and the differential $b:=\sum_{i=0}^{n+1}(-1)^i\Delta_i$.
\item The norm $N:=\sum_{i=0}^nT^i:C^{\ot n+1}\to C^{\ot n+1}$
\end{itemize}
In the same way as in the non-differential graded case, we have
the following Lemma and Theorem:
\begin{lem}
With the above notations:
$$T^{n+1}_{C^{\ot n+1}}=id_{C^{\ot n+1}}\ \ ; \ \ \ Nb'=bN\ $$
$$(1-T)N=0=N(1-T) \ ; \ (1-T)b=b'(1-T)$$
\end{lem}
\dem They are formal consequence of the relations
$$\left\{
\begin{array}{cc}
T\Delta_1=-\Delta_{i-1}T& i=1,\dots, n\\
T\Delta_0=(-1)^{n+1}\Delta_{n+1}&
\end{array}
\right.$$
We also have that the signs on $T$ are defined in such a way that
$d$ commutes with $T$, so $d$ commutes with $(1-T)$ and $N$. Then there
is a well-defined double complex, denoted by
$\C^{**}(C)$:
{\small$$\diagram
\dots&
\C_{\Hoch}(C)[-4]\lto_{1-T\ \ } &
\wh{C}(C)[-3]\lto_{N}&
\C_{\Hoch}(C)[-2]\lto_{1-T} &
\wh{C}(C)[-1]\lto_{N}&
\C_{\Hoch}(C)\lto_{1-T} 
\enddiagram$$}
where $\C_{\Hoch}(C)=Tot(\oplus_{n\geq 0}C^{\ot n+1},b,d)$ and 
$\wh{C}(C)=Tot(\oplus_{n\geq 0}C^{\ot n+1},b',d)$. We define, for a
differential graded coalgebra $C$, $\Hoch^*(C):=H^*(\C_{\Hoch}(C))$,
and $HC^*(C):=H^*(Tot(\C^{**}(C)))$.

\bigskip

\rem The difference between $\wh{C}(C)$ and the standard
resolution $C(C)$ is that $\wh{C}(C)=
Tot(\oplus_{n\geq 0}C^{\ot n+1},b',d)$
while $C(C)=Tot(\oplus_{n\geq 1}C^{\ot n+1},b',d)$. In fact
$\wh{C}(C)$ can be identified with the mapping cone
of $\rho_C:C\to C(C)$ which is a quasi-isomorphism, then
$\wh{C}(C)$ is an acyclic complex.

\begin{teo}\label{sbi}
Let $C$ be a differential graded coalgebra, then
there is an SBI-type long exact sequence
$$\diagram 
\dots \rto &HC^n(C)\rto^S &HC^{n+2}(C)\rto^I &
\Hoch^{n+2}(C)\rto^B &HC^{n+1}(C)\rto &\dots 
\enddiagram$$
\end{teo}
\dem we proceed in the same way as in the non differential graded case,
noticing the 2-periodicity in the definition of
$\C^{**}(C)$ we obtain a short exact sequence
of double complexes
$$0\to \C^{**}(C)[-2]\to \C^{**}(C)\to Co(1-T)[1]\to 0$$
and hence a long exact sequence on cohomology.

Since
$\wh{C}(C)$ is acyclic, it follows that $Co(1-T)[1]$
is quasi-isomorphic to $\C_{Hoch}(C)$. 

\begin{lem}
Let $C$ be a positively graded differential 
coalgebra, then the standard resolution $C(C)$ is a closed object 
in $\Chain(C^e)$.
\end{lem}
\dem it follows the lines of the proof that $C(M)$ is a closed object
in $\Chain(C)$ for all $M\in\Chain^+(C)$. Since $C$ is positively graded,
then $C(C)=\stackunder{\lo p}{\lim}C_{\leq p}(C)$ where
$C_{\leq p}(C)=Tot(\oplus_{j=1}^pC^{\ot j+1},b',d)$ is the quotient of
$C(C)$ by the subcomplex $Tot(\oplus_{j>p}C^{\ot j+1},b',d)$, and
we have short exact sequences
$$0\to (C^{\ot p+2},d)\to C_{\leq p+1}(C)\to C_{\leq p}(C)\to 0$$
which clearly split as sequences of graded $C$-bicomodules.
By the isomorphism $(C^{\ot p+2},d)\cong C^e\ot C^{\ot p}$, it follows
that $(C^{\ot p+2},d)$ is closed in $\Chain(C^e)$ and
$C_{\leq 1}(C)=C\ot C=C^e$ which is also $C^e$-closed, then inductively
$C_{\leq p}(C)$ is a closed $C^e$-object for all $p\geq 1$.
Now using again the above exact sequence together with 
Lemma \ref{tri} we can conclude that $\stackunder{\lo p}{\lim}C_p(C)$
is closed in $\Chain(C^e)$.

\begin{coro}
Let $C$ be a positively graded differential coalgebra, then
$\C_{\Hoch}(C)=C\Box_{C^e}^RC$.
\end{coro}
\dem in the same way as in the non-differential graded case, it
is easy to see that $\C_{\Hoch}(C)\cong C(C)\Box_{C^e}C$, the
equality $C(C)\Box_{C^e}C=C\Box_{C^e}^RC$ comes from the facts
that $C$ is quasi-isomorphic to $C(C)$ (as $C$-bicomodule)
and that $C(C)$ is a closed
object in $\Chain(C^e)$.

\begin{defi}
Let $C$ be a differential graded coalgebra and $M\in \Chain(C^e)$,
we define
$$\Hoch^n(M,C):=H^n(C(C)\Box_{C^e}M)$$
$$H^n(M,C):=\Hom_{\D(C^e)}(M,C[n])$$
\end{defi}

\begin{teo}\label{morita}
Let $C$ and $D$ be two positively graded coalgebras such that there
exist $P\in\Chain^+(C\ot D^{op})$ and
 $Q\in\Chain^+(D\ot C^{op})$ with
$P\Box_D^RQ\cong C$ and $Q\Box_C^RP\cong D$ (isomorphisms in
$\D(C^e)$ and $\D(D^e)$ respectively). Then, for all $M\in\Chain^+(C)$
$$H^*(M,C)\cong \Hoch(Q\Box_C^RM\Box_C^RP,D)$$
$$H^*(M,C)\cong H^*(Q\Box_C^RM\Box_C^RP,D)$$
In particular $\Hoch^*(C)\cong \Hoch^*(D)$ and $H^*(C)\cong H^*(D)$.
\end{teo}
\dem It is clear that
$Q\Box_C^R-\Box_C^RP:\D(C^e)\to \D(D^e)$ is an equivalence of
triangulated categories, then
$$H^n(M,C)=\Hom_{\D(C^e)}(M,C[n])\cong
\Hom_{\D(D^e)}(Q\Box_C^RM\Box_C^RP,Q\Box_C^RC[n]\Box_C^RP)\cong$$
$$\cong \Hom_{\D(D^e)}(Q\Box_C^RM\Box_C^RP,D[n])=
H^n(Q\Box_C^RM\Box_C^RP,D)$$
For the other cohomology theory,
$$\Hoch^*(Q\Box_C^RM\Box_C^RP,D)=
H^*\left((Q\Box_C^RM\Box_C^RP)\Box_{D^e}^RD\right)\cong$$
$$\cong H^*\left(M\Box_{C^e}^R(P\Box_{D}^RD\Box_D^RQ)\right)\cong
H^*(M\Box_{C^e}^RC)=\Hoch^*(M,C)$$

\bigskip

\rem $H^*(C)$ is a graded algebra, with multiplication given by
composition in $\D(C^e)$. In the situation of the above theorem,
since the isomorphism $H^*(C)\cong H^*(D)$ is given by a functor,
it follows that the isomorphism is not only an isomorphism
of $k$-vector spaces
but also of graded $k$-algebras.


\bigskip

\exs 1. Let $C$ and $D$ be two usual coalgebras such that they
admit a cotilting bicomodule ${}_DT_C$, then $\Hoch^*(C)\cong
\Hoch^*(D)$ and $H^*(C)\cong H^*(D)$.

\medskip

\noindent 2. Let $f:C\to D$ be a quasi-isomorphism of positively 
graded differential coalgebras, then $\Hoch^*(M,C)\cong
\Hoch^*({}_fM_f,D)$ and $H^*(M,C)\cong H^*({}_fM_f,D)$ (where
${}_fM_f$ is $M$, but viewed as $D$-bicomodule via $f$).

\bigskip

In the case of the second example, we also have that the 
quasi-isomorphism $C\Box_{C^e}^RC\cong D\Box_{D^e}^RD$ is
induced by $f$, we
will see next that this implies the invariance of $HC^*$:

\begin{prop}\label{propf}
Let $C$ and $D$ be two differential graded
coalgebras and $f:C\to D$ be a 
differential graded and coassociative map
such that
$f_*:\C_{\Hoch}(C)\to\C_{\Hoch}(D)$ is a quasi-isomorphism,
then $f$ induces an isomorphism $HC^*(C)\cong HC^*(D)$.
\end{prop}
\dem we begin by noticing that the complex $\C^{**}(C)$
(also $C(C)$ and $\C_{\Hoch}(C)$) is defined using the differential graded
structure and the comultiplication, so $f$ induces a  morphism
of complexes $\C^{**}(C)\to\C^{**}(D)$ which is natural
with respect to the SBI-long exact sequence, and we obtain
a morphism of long exact sequences:
$$\diagram 
\dots \rto &HC^n(C)\rto^S\dto^{f} &HC^{n+2}(C)\dto^{f}\rto^I &
\Hoch^{n+2}(C)\rto^B \dto^{\wt{\ }}&HC^{n+1}(C)\dto^f\rto &\dots \\
\dots \rto &HC^n(D)\rto^S &HC^{n+2}(D)\rto^I &
\Hoch^{n+2}(D)\rto^B &HC^{n+1}(D)\rto &\dots 
\enddiagram$$
So the proposition is proved by induction, using that
$HC^0=\Hoch^0$ and the five Lemma for the inductive step.

\bigskip

As a corollary we have that if $f:C\to D$ is
a quasi-isomorphism of positively graded differential coalgebras
then $HC^*(C)\cong HC^*(D)$.  We remark that if $C$ or $D$ have
some nonzero component in negative degrees, we do not know
if the standard resolution is a closed object, we cannot
interpret $C(C)\Box_{C^e}C$ as $C\Box_{C^e}^RC$  and so we 
cannot apply Proposition \ref{propquis}.
In fact we do not know even how to define 
$-\Box_{C}^R-$ in general when the coalgebra $C$ is not positively graded
because we do not know if there exist ``enough" closed objects.

On the other hand, it would be interesting to know if 
$HC^*$ is invariant under cotilting equivalence. This 
would be implied if one shows that the coassociative maps 
$e_D(T\oplus D)\to e_D(T)=C$ and $e_D(T\oplus C)\to e_D(D)=D$ 
(where ${}_DT_C$ is a cotilting bicomodule) induce quasi-isomorphisms 
in $\C_{\Hoch}$.  We know (see \cite{hc}) that this is the case 
when $C$ and $D$ are Morita - Takeuchi equivalent.

\end{document}